\documentclass[12pt,reqno]{amsart}
\usepackage{mathrsfs}
\usepackage{amsgen}
\usepackage{amscd}

\usepackage{amsmath}
\usepackage{latexsym}
\usepackage{amsfonts}
\usepackage{amssymb}
\usepackage{amsthm}
\usepackage{graphicx,epsfig}
\usepackage{color}
\usepackage{amsthm,amssymb}
\usepackage{graphicx}
\usepackage{enumerate}
\usepackage{amssymb,amsmath,graphicx,amsfonts,euscript}
\usepackage{color}
\usepackage{mathrsfs}
\usepackage{ulem}
\usepackage[margin=3cm, a4paper]{geometry}
\DeclareGraphicsExtensions{.eps}
\usepackage{amssymb,amsmath,graphicx,amsfonts,euscript,mathrsfs}
\usepackage{hyperref}

\DeclareGraphicsExtensions{.pdf}

\usepackage{epstopdf}
\newtheorem{theorem}{Theorem}[section]
\newtheorem{lemma}[theorem]{Lemma}
\newtheorem{proposition}[theorem]{Proposition}

\theoremstyle{definition}

\theoremstyle{remark}

\def\R{\mathbb{R}}
\def\Z{\mathbb{Z}}
\def\T{\mathbb{T}}
\def\C{\mathbb{C}}

\newcommand{\fe}{\mathrm{e}}

\newcommand{\bR}{{\mathbb R}}

\newcommand{\bT}{{\mathbb T}}

\numberwithin{equation}{section}

\allowdisplaybreaks

\begin{document}
\makeatletter
\@namedef{subjclassname@2020}{\textup{2020} Mathematics Subject Classification}
\makeatother

\title[Fully discrete integrator for  NLS equation]{A fully discrete low-regularity integrator for the nonlinear Schr\"odinger equation}

\author[A. Ostermann]{Alexander Ostermann}
\address{\hspace*{-12pt}A.~Ostermann: Department of Mathematics, University of Innsbruck, Technikerstrasse 13, 6020 Innsbruck, Austria}
\email{alexander.ostermann@uibk.ac.at}

\author[F. Yao]{Fangyan Yao}
\address{\hspace*{-12pt}F.~Yao: School of Mathematical Sciences,
South China University of Technology,
 Guangzhou, Guangdong 510640, P. R. China}
\email{fangyanyao@outlook.com}

\subjclass[2020]{Primary 65M12, 65M15, 35Q55. }

\keywords{  Low regularity, nonlinear Schr\"{o}dinger equation, fully discrete, fast Fourier transform}

\maketitle

\begin{abstract}\noindent
For the solution of the one dimensional cubic nonlinear Schr\"odinger equation on the torus, we propose and analyse a fully discrete low-regularity integrator. The considered scheme is explicit. Its implementation relies on  the fast Fourier transform with a complexity of $\mathcal{O}(N\log N)$ operations per time step,
 where $N$ denotes the degrees of freedom in the spatial discretisation.
  We prove that the new scheme provides an
$\mathcal{O}(\tau^{\frac32\gamma-\frac12-\varepsilon}+N^{-\gamma})$ error bound  in $L^2$ for any initial data in $H^\gamma$,  $\frac12<\gamma\leq 1$, where $\tau$  denotes the temporal step size.
Numerical examples illustrate this convergence behavior.
\end{abstract}

\section{Introduction}\label{sec:introduction}
In this paper, we analyze low-regularity integrators for
 the   cubic nonlinear Schr\"odinger equation  (NLS):
\begin{equation}\label{model}
	\left\{\begin{aligned}
		& i\partial_tu(t,x)+\partial_{xx} u(t,x)
		=|u(t,x)|^2u(t,x),\\
		&u(0,x)=u^0(x).
	\end{aligned}\right.
\end{equation}
We consider this problem on the one dimensional  torus $\T=(0,2\pi)$.
The function $u:\bR^{+}\times\bT\to\C$ is the unknown
and $u^0\in H^\gamma(\bT)$, $\gamma\geq0$  is the given initial data.
The NLS equation \eqref{model} is globally well-posed in $H^\gamma(\T)$, $\gamma\geq0$; see, e.g., \cite{Bo}.

The construction of efficient numerical schemes for dispersive equations has been the subject of much research work.
In particular, for the NLS equation,
  substantial research has been undertaken in numerical analysis.
For smooth solutions,  many classical numerical methods have been analyzed for  numerical discretisation in space and time,
 for example, finite difference methods \cite{sanz}, operator splitting \cite{besse,jahnke,Lubich}, and exponential integrators \cite{faou}.
 These methods have all their own characteristics.
 However, they generally require relatively high regularity
 of the solution.
 For example, initial data in $H^{\gamma+2}$ are required to obtain
 first-order convergence in $H^\gamma$ for the NLS equation in \cite{Lubich}, and   initial data in $H^{\gamma+4}$ for second-order convergence.

 As mentioned above, if the exact solution is smooth enough,
one can rely on   classical numerical schemes.
In practical application, however, nonsmooth initial data
are encountered as well. A typical example is applications in nonlinear optics, where initial data can be corrupted with noise.
Therefore,  recent attention has been focused on  lower regularity requirements.
In order to achieve convergence with the lowest possible regularity
of the initial data,
so-called low-regularity integrator were proposed recently.

For the NLS equation, Ostermann and Schratz \cite{lownls} achieved  first-order convergence
in $H^\gamma(\T^d)$, $d\geq1$ with initial data $H^{\gamma+1}(\T^d)$
by introducing  a new exponential-type numerical scheme
under the assumption $\gamma>\frac{d}{2}$.
Furthermore, the authors obtained the convergence in  $H^\gamma(\T)$ with initial data $H^{\gamma}(\T)$
for the quadratic NLS equation in  one space dimension.
Then Wu and Yao \cite{w} constructed a new first-order scheme
for the cubic NLS equation in one space dimension and obtained  first-order convergence in
$H^\gamma(\T)$ with initial data $H^{\gamma}(\T)$ for $\gamma>\frac32$.
A second-order scheme was proposed by Kn\"{o}ller, Ostermann and Schratz in \cite{lownls2}. In one space dimension, this scheme  requires two additional derivatives of the solution; in higher dimensions, three additional derivatives  are necessary.
Later, Ostermann, Rousset and Schratz \cite{ostermann,ostermann2020}
proved  convergence in  $L^2$ for  initial data
in $H^s(\R^d)$ and $H^s(\T)$ respectively, $0<s\leq1$ and obtained
  fractional orders of convergence in a frame work of discrete Bourgain spaces.
 All results discussed above concern time integration only.
Recently, Li and  Wu \cite{liwu} considered a fully discrete low-regularity integrator in one space dimension
for the NLS equation and got  first-order convergence (up to a logarithmic factor)
 in both time and space in $L^2(\T)$ for $H^1(\T)$ initial data.

The purpose of this article is to construct a fully discrete low-regularity integrator and to prove convergence
for initial data $u^0\in H^s(\T)$, $\frac12<s\leq 1$.
For the spatial discretisation, we use a Fourier ansatz with frequency truncation.
For the temporal discretisation,
we employ a careful analysis of
the nonlinear dynamics in phase space.
In addition, we use harmonic analysis technique in the proof of
some technical lemmas.

The rest of this paper is structured as follows.
 We present the fully discrete low-regularity integrator in section 2. We state the main convergence result
 and explain the employed notations.
In section 3, we derive the considered  scheme and give some technical lemmas,
which will be used in the convergence proof.
In section 4, we show the $H^\gamma$ error bound before the poof of the theorem.
Finally we prove stability and  the error bound  of the fully discrete scheme.
In sections 5 and 6, we report numerical experiments that illustrate
our theoretical analysis and we draw some conclusions.

\section{Notations and main result}

\subsection{Some notations}\label{subsec1}
 We denote by $\langle \cdot,\cdot\rangle$  the $L^2$ inner product
 on $\T$, that is
$$
\langle f,g \rangle=\int_\T f(x) \overline{g(x)}\,dx, \qquad f, \ g\in L^2(\T;\C).
$$
The Fourier transform $(\hat{f}_k)_{k\in \Z}$ of a function $f$: $\T\rightarrow\C$ is defined by
$$
\hat{f}_k=\frac{1}{2 \pi}\int_{\T} e^{- i   kx }f( x )\,d x.
$$
The Fourier inversion formula is given by
$$
f( x )=\sum_{k\in \Z} \hat{f}_k e^{ i  kx } .
$$
We recall the following properties:
\begin{eqnarray*}
 &\|f\|_{L^2(\T)}^2
 = 2 \pi\sum\limits_{k\in \Z}\big|\hat{f_k}\big|^2,\quad f\in L^2(\T); \\
 & \widehat{(fg)}_k=\sum\limits_{k=k_1+k_2}
  \hat{f}_{k_1}\hat{g}_{k_2},\quad f, \ g\in L^2(\T).
\end{eqnarray*}
For our convergence analysis, we equip the  Sobolev space $H^s(\T)$, $s\geq0$ with the norm
$$
\big\|f\big\|_{H^s(\T)}^2=
\big\|J^sf\big\|_{L^2(\T)}^2
=2 \pi\sum_{k\in \Z}(1+ k ^2)^s |\hat{f}_k|^2, \qquad J^s=(1-\partial_{xx})^\frac s2.
$$
This norm is equivalent with the standard norm of $H^s(\T)$.

Further, we denote by $\partial_x^{-1}$ the operator defined in Fourier space as
\begin{equation}\label{def:px-1}
(\widehat{\partial_x^{-1}f})_k
=\Bigg\{ \aligned
    &(i k )^{-1}\hat{f}_k \quad &\mbox{if  }  k \ne 0,\\
    &0 \quad &\mbox{if  }  k = 0.
   \endaligned
\end{equation}

For convenience, we introduce some additional notations.
First, we define  the zero-mode operator by
$$
 P_0f=\hat f_0=\frac{1}{2 \pi}\int_\T f(x)\,dx.
$$
Furthermore, for any positive integer $N$, we define the projection operators $ P_N$
and $ P_{>N}$ by
\begin{align}\label{truncation}
(\widehat{ P_N f})_k
= \left\{\begin{aligned}
& \hat f_k
 &\quad \mbox{if  }  |k|\leq N,\\
 &0 &\quad \mbox{if  }  |k|>N;
 \end{aligned}\right.
 \qquad
( \widehat{ P_{>N} f})_k
= \left\{\begin{aligned}
& \hat f_k
 &\quad \mbox{if  }  |k|> N,\\
 &0 &\quad \mbox{if  }  |k|\leq N.
 \end{aligned}\right.
\end{align}
Let $S_N$ be the space consisting of all functions  $f\in L^2(\T)$ such that $\hat f_k=0$ for $|k|>N$.
Then, for $f$ and $g\in S_N$, the cost of computing the Fourier coefficients of $ P_N(fg)\in S_N$ is $\mathcal{O}(N\log N)$.
This can be seen as follows.

Let $I_{2N}$ be the $(4N+1)$-point trigonometric interpolation operator
\begin{align}\label{in}
I_{2N}f(x)=\sum_{k=-2N}^{2N}\fe^{ikx}\tilde{f}_k,
\end{align}
where
$$
\tilde{f}_k=\frac{1}{4N+1}\sum_{k=-2N}^{2N}\fe^{-ikx_n}f(x_n),
\quad x_n=\frac{2 \pi n}{4N+1},\quad n=-2N,\cdots,2N.
$$
Consequently, if the Fourier coefficient $\hat f_k$ of the function $f$ satisfies  $\hat f_k=0$ for $|k|>2N$,
then $I_{2N} f=f$ and $\tilde{f}_k=\hat f_k$.
Furthermore,  if $f$ and $g\in S_N$,
then $\widehat{(fg)}_k=0$ for  $|k|>2N$.
Hence  we get $fg=I_{2N}(fg)$ from the definition of $I_{2N}$.
Further, the cost of computing the Fourier coefficients of $ P_NI_{2N}(fg)\in S_N$ is $\mathcal{O}(N\log N)$.

\subsection{Numerical method and main result}
Let $\tau$ denote the temporal step size and
$t_n=n\tau$ the corresponding  sequence of
discretisation points in the time interval $[0,T]$.
In section 3.1 below, we construct the fully discrete low-regularity integrator for  equation \eqref{model} as
\begin{align}\label{numerical}
u^{n+1}_{\tau,N}=\Phi_{\tau,N}(u^n_{\tau,N})
\end{align}
with $u^0_{\tau,N}= P_N I_{2N}u_0$.
The numerical flow $\Phi_{\tau,N}$ maps a function $f\in S_N$ to
 $\Phi_{\tau,N}(f)\in S_N$. It is defined in the following way:
\begin{align}\label{psi}
\Phi_{\tau,N}(f)
&=\fe^{i\tau\partial_x^2}f
+\frac12 \partial_x^{-1} P_N
          \Big[ \big(\fe^{-i\tau\partial_x^2} \partial_x^{-1}\bar f\big)
              \cdot \fe^{i\tau\partial_x^2} P_N
                \big( f^2\big)  \Big]
 -\frac12\fe^{i\tau\partial_x^2} \partial_x^{-1} P_N
          \Big[ \partial_x^{-1} \bar f\cdot  P_N\big( f^2\big)   \Big]\notag\\
&\qquad -i\tau P_0\Big[ \bar f
              \cdot  P_N
                \big( f^2\big)\Big]
-i\tau\Big[
\fe^{i\tau\partial_x^2} P_N
\big( f^2\big)
- P_0(f^2)\Big]
\cdot  P_0(\bar f)\notag\\
 &\qquad+\frac12 \fe^{i\tau\partial_x^2}  P_N
 \Big[\bar f
  \cdot \fe^{-i\tau\partial_x^2}
  P_N \big(\fe^{i\tau\partial_x^2}
  \partial_x^{-1} f\big)^2\Big]
-\frac12 \fe^{i\tau\partial_x^2}  P_N
 \Big[\bar f\cdot
  P_N \big(\partial_x^{-1} f\big)^2\Big]\notag\\
  &\qquad+i\tau\fe^{i\tau\partial_x^2}  P_N
 \Big[\bar f
  \cdot
  P_N \big( f^2\big)\Big]
 -2i\tau\fe^{i\tau\partial_x^2}  P_N
 \Big[\bar f\cdot P_N f\Big]
 \cdot  P_0(f)\notag\\
 &\qquad+i\tau\fe^{i\tau\partial_x^2}  P_N {\bar f}
 \cdot \big( P_0(f))^2.
\end{align}
From the discussion in section \ref{subsec1},  we infer that the initial data $u^0_{\tau,N}= P_N I_{2N}u_0$ can be computed with FFT with computational
cost of $\mathcal{O}(N\log N)$. Furthermore, the terms $\fe^{i\tau\partial_x^2}f$
and $P_N(fg)$ also can be computed with FFT for  given functions $f,g\in S_N$.
We note that \eqref{psi} only consists of the above expressions.  Therefore, \eqref{psi} can be computed with FFT.

Now, we state the main result of this paper. We show  convergence of the fully discrete low-regularity integrator
given in \eqref{numerical}.
\begin{theorem}\label{result1}
Let $u^n_{\tau,N}$ be the numerical solution \eqref{numerical} of   equation \eqref{model}
up to some fixed time $T>0$. Under the assumption that $u^0\in H^{\gamma}(\bT)$ for some $\frac12<\gamma\leq1$,
  for arbitrary given  $\varepsilon>0$, there exist constants $\tau_0$, $C>0$ such that for  any
  step size $0<\tau\leq\tau_0$  and all $0\leq n\tau\leq T$
\begin{equation}
  \|u(t_n,\cdot)-u^n_{\tau,N}\|_{L^2}\leq C\tau^{\frac32\gamma-\frac12-\varepsilon}
+C N^{-\gamma},
\end{equation}
 where the constant $\tau_0$  depends only on $T$ and $\|u\|_{L^\infty((0,T);H^{\gamma})}$,
 and the constant $C$  depends only on $T$, $\|u\|_{L^\infty((0,T);H^{\gamma})}$ and $\varepsilon$.
\end{theorem}

We write $A\lesssim B$ or $B\gtrsim A$ to express that $A\leq CB$ for some positive  constant $C$ which may
be different at each constant but is independent of $\tau$, $N$ and $n$. Further,  we write $A\sim B$ for
$A\lesssim B\lesssim A$. We write $\mathcal{O}(Y)$ to denote a quantity $X$ that satisfied $|X|\lesssim |Y|$.

Henceforth, we denote by $ T_m(M; v)$ the class of  functions $f\in L^2(\T)$ such that
\begin{align}\label{def:TmM}
|\hat f_k|\lesssim
\sum\limits_{ k = k _1+\cdots+ k _m}
| M( k, k _1,\cdots, k _m)|\>|\hat{v}_{ k _1}| \cdots |\hat v_{ k _m}|,
\quad {\rm for\ all}\quad k,
\end{align}
where $v_k$ is the $k$th Fourier coefficient of $v$.

\section{The construction of the scheme and some technical lemmas}

In this section, we construct the fully discrete low-regularity
exponential integrator by frequency truncation and harmonic analysis techniques. Then, we state some lemmas that will be used  frequently
in section 4.

\subsection{The construction of the scheme}
It is known that if $u^0\in H^\gamma(\T)$, $\gamma\geq0$, then the NLS equation \eqref{model} has a unique solution $u\in C([0,T];H^\gamma(\T))$; see \cite{Bo}.
Recalling  Duhamel's formula, we write
\begin{align*}
u(t_{n+1})=\fe^{i\tau\partial_x^2}u(t_n)
 -i\int_0^\tau \fe^{i\left(t_{n+1}-(t_n+s)\right)\partial_x^2}
   \big[|u(t_n+s)|^2u(t_n+s)\big]\,ds.
\end{align*}
With the twisted variable $v(t)=\fe^{-it\partial_x^2}u(t)$, the above formula becomes
\begin{align}\label{solution}
v(t_{n+1})=v(t_n)-i\int_0^\tau \fe^{-i(t_n+s)\partial_x^2}
    \big[|\fe^{i(t_n+s)\partial_x^2}v(t_n+s)|^2\,\fe^{i(t_n+s)\partial_x^2}v(t_n+s)\big]\,ds.
\end{align}
Applying the Fourier transform, \eqref{solution} can be expressed as
\begin{align*}
\hat v_{ k }(t_{n+1})
  =\hat v_{ k }(t_n)
     -i\int_0^\tau\sum\limits_{ k = k _1+ k _2+ k _3}\fe^{i(t_n+s)\phi}
        \>\widehat{ \bar v}_{ k _1}(t_n+s)\hat v_{ k _2}(t_n+s)\hat v_{ k _3}(t_n+s)\,ds,
\end{align*}
where $\hat v_{ k }(t)$ denotes the $k$th Fourier coefficient of $v(t)$.  We also use the phase function
\begin{align}\label{def-1}
\phi( k, k _1, k _2, k _3)= k ^2+ k _1^2- k _2^2- k _3^2.
\end{align}
In order to obtain a first-order approximation,
by \eqref{solution}, we have for any $s\in [0,\tau]$,
\begin{align}
v(t_n+s)-v(t_n)\in  \tau T_3(1;v).
\end{align}
This implies that
\begin{align}\label{v-1st-phi}
\hat v_{ k }(t_{n+1})
=\hat v_{ k }(t_n)+\hat{I}_{1,k} +\hat{\mathcal{R}}_{1,k}(v),
\end{align}
where
\begin{align*}
\hat{I}_{1,k}=-i\sum\limits_{ k = k _1+ k _2+ k _3}\int_0^\tau \fe^{i(t_n+s)\phi}\,ds
  \>\widehat{ \bar v}_{ k _1}\hat v_{ k _2}\hat v_{ k _3}\quad
{\text{and}}\quad \mathcal{R}_1(v)\in\tau^2 \> T_5(1;v).
\end{align*}
Henceforth, we denote   $\hat v_{ k }(t_n)$ by $\hat v_{ k }$ and
  the  $k$th Fourier coefficient of $\mathcal{R}_j(v)$ for $j\geq 1$ by $\hat{\mathcal{R}}_{j,k}(v)$ for short.

For further approximation, we consider a decomposition into
low and high frequencies. In particular, we consider
 the following two cases: $|k|\leq N$ and $|k|> N$.

Case 1: $|k|\leq N$.  We consider only the first  term $\hat{I}_{1,k}$ in \eqref{v-1st-phi} and
 truncate $\hat{I}_{1,k}$ to the frequency domain $|k_2+k_3|\leq N$,
\begin{equation}\label{for4}
  \hat{I}_{1,k} =\mathcal{K}_k(v)+\hat{\mathcal{R}}_{2,k}(v),
\end{equation}
where $\mathcal{K}_k(v)$ is defined by
$$
\mathcal{K}_k(v)=-i\sum\limits_{\substack{k=k_1+k_2+k_3\\|k_2+k_3|\leq N}}\int_0^\tau \fe^{i(t_n+s)\phi}\,ds
       \>\widehat{ \bar v}_{k_1}\hat v_{k_2}\hat v_{k_3}.
$$
 The remainder
$\hat{\mathcal{R}}_{2,k}(v)$ is given by
\begin{equation}
\hat{\mathcal{R}}_{2,k}(v)=
 -i\sum\limits_{{\substack{k=k_1+k_2+k_3\\   |k_2+k_3|> N}}}\int_0^\tau \fe^{i(t_n+s)\phi}\,ds
       \>\widehat{ \bar v}_{k_1}\hat v_{k_2}\hat v_{k_3}.
\end{equation}
Furthermore, from the definition of $ T_m(M;v)$ in \eqref{def:TmM},
the function $\mathcal{R}_2$  satisfies
$$
\mathcal{R}_{2}(v)\in\tau \> T_3(1_{|k_2+k_3|> N};v).
$$
Next we consider the term $\mathcal{K}_k(v)$.
Note that if $k=k_1+k_2+k_3$, then from \eqref{def-1}, the following equality holds
\begin{align*}
\phi( k, k _1, k _2, k _3)
=2kk_1+2k_2k_3.
\end{align*}
In order to get a first-order scheme,
we need to find an appropriate approximation to the exponential $\fe^{is\phi}$.
Using the formula
\begin{align*}
\fe^{is\phi}=
\fe^{2iskk_1}\fe^{2isk_2k_3}
 =\fe^{2iskk_1}+\big(\fe^{2isk_2k_3}-1\big)
 +\big(\fe^{2iskk_1}-1\big)\big(\fe^{2isk_2k_3}-1\big),
\end{align*}
we decompose $\mathcal{K}_k(v)$ into two terms
\begin{subequations}
\begin{align}
\mathcal{K}_k(v)
=-i\sum\limits_{{\substack{k=k_1+k_2+k_3\\   |k_2+k_3|\leq N}}}
\int_0^\tau \fe^{it_n\phi}\Big[\fe^{2iskk_1}+\big(\fe^{2isk_2k_3}-1\big)\Big]\,ds
       \>\widehat{ \bar v}_{k_1}\hat v_{k_2}\hat v_{k_3}
+\hat{\mathcal{R}}_{3,k}(v)\notag,
\end{align}
\end{subequations}
where $\hat{\mathcal{R}}_{3,k}(v)$ is defined by
\begin{equation}\label{def:r3}
\hat{\mathcal{R}}_{3,k}(v)
= -i\sum\limits_{{\substack{k=k_1+k_2+k_3\\   |k_2+k_3|\leq N}}}\int_0^\tau \fe^{it_n\phi}
 \big(\fe^{2iskk_1}-1\big)\big(\fe^{2isk_2k_3}-1\big)\,ds
       \>\hat{ \bar v}_{k_1}\hat v_{k_2}\hat v_{k_3}.
\end{equation}

We find that the integral in $\mathcal{K}_k(v)$ can be computed easily.
 Integrating with respect to $s$, we have for any $|k|\leq N$,
\begin{align}\label{kapa}
&\mathcal{K}_k(v)\notag\\
  =&  -i\sum\limits_{{\substack{k=k_1+k_2+k_3\\   |k_2+k_3|\leq N
     \\|k|\neq 0,|k_1|\neq 0}}}\frac{1}{2ikk_1}\fe^{it_n\phi}
\big(\fe^{2i\tau kk_1}-1\big)
       \>\widehat{ \bar v}_{k_1}\hat v_{k_2}\hat v_{k_3}
-i\tau\sum\limits_{{\substack{0=k_1+k_2+k_3\\   |k_2+k_3|\leq N}}}\fe^{it_n(k_1^2-k_2^2-k_3^2)}
       \>\widehat{ \bar v}_{k_1}\hat v_{k_2}\hat v_{k_3}\notag\\
&-i\tau\sum\limits_{{\substack{k=k_2+k_3\\   |k_2+k_3|\leq N}}}\fe^{it_n(k^2-k_2^2-k_3^2)}
       \>\hat v_{k_2}\hat v_{k_3}\widehat{ \bar v}_{0}
+i\tau\sum\limits_{0=k_2+k_3}
 \fe^{it_n(-k_2^2-k_3^2)}
       \>\hat v_{k_2}\hat v_{k_3}\widehat{ \bar v}_{0}\notag\\
 &-i\sum\limits_{{\substack{k=k_1+k_2+k_3\\   |k_2+k_3|\leq N\\|k_2|\neq 0,|k_3|\neq 0}}}\frac{1}{2ik_2k_3}\fe^{it_n\phi}
\big(\fe^{2i\tau k_2k_3}-1\big)
       \>\widehat{ \bar v}_{k_1}\hat v_{k_2}\hat v_{k_3}
+i\tau\sum\limits_{\substack{k=k_1+k_2+k_3\\   |k_2+k_3|\leq N}}\fe^{it_n\phi}
       \>\widehat{ \bar v}_{k_1}\hat v_{k_2}\hat v_{k_3}\notag\\
&-2i\tau\sum\limits_{\substack{k=k_1+k_2\\   |k_2|\leq N}}\fe^{it_n(k^2+k_1^2-k_2^2)}
       \>\widehat{ \bar v}_{k_1}\hat v_{k_2}\hat v_{0}
+i\tau\fe^{2it_n k^2}\>\widehat{ \bar v}_{k}
\big(\hat v_{0}\big)^2
+\hat{\mathcal{R}}_{3,k}(v).
\end{align}
For specific details of the above formula, we  refer to  the literature \cite{lownls2, liwu}.

Case 2: $|k|> N$.
Let $\mathcal{R}_4(v)$ be the function with Fourier coefficients
\begin{align*}
\hat{\mathcal{R}}_{4,k}(v)
=-i\sum\limits_{ \substack{k = k _1+ k _2+ k _3\\|k|>N}}\int_0^\tau \fe^{i(t_n+s)\phi}\,ds
       \>\widehat{ \bar v}_{ k _1}\hat v_{ k _2}\hat v_{ k _3}.
\end{align*}
Then, for $|k|>N$, we get
\begin{equation}\label{high}
  \hat{I}_{1,k} =\hat{\mathcal{R}}_{4,k}(v).
\end{equation}
Furthermore, from the definition of $ T_m$, we find
$$
\mathcal{R}_{4}(v)\in\tau \> T_3(1_{|k|> N};v).
$$
Hence, putting together  \eqref{v-1st-phi}, \eqref{for4},  \eqref{kapa} and \eqref{high} yields
\begin{align}\label{vtn+1-app}
v(t_{n+1})=
&\Phi_{\tau,N}^n\big(v(t_n)\big)+\mathcal{R}_1(v)
+\mathcal{R}_2(v)+\mathcal{R}_3(v)+\mathcal{R}_4(v),
\end{align}
where $\Phi_{\tau,N}^n$ is given by
\begin{align}\label{Phi-def}
\Phi_{\tau,N}^n(f)
&=f+\frac12\fe^{-it_{n+1}\partial_x^2} \partial_x^{-1} P_N
          \Big[ \big(\fe^{-it_{n+1}\partial_x^2} \partial_x^{-1} \bar f\big)
              \cdot \fe^{i\tau\partial_x^2} P_N
                \big(\fe^{it_n\partial_x^2} f\big)^2   \Big]\notag\\
&\quad -\frac12\fe^{-it_n\partial_x^2} \partial_x^{-1} P_N
          \Big[ \big(\fe^{-it_n\partial_x^2} \partial_x^{-1}\bar f\big)
              \cdot  P_N
                \big(\fe^{it_n\partial_x^2} f\big)^2   \Big]\notag\\
&\quad -i\tau P_0\Big[\big(\fe^{-it_n\partial_x^2} \bar f\big)
              \cdot  P_N
                \big(\fe^{it_n\partial_x^2} f\big)^2\Big]\notag\\
&\quad -i\tau\fe^{-it_n\partial_x^2} P_N
\big(\fe^{it_n\partial_x^2} f\big)^2
\cdot  P_0(\bar f)
+i\tau P_0
\big(\fe^{it_n\partial_x^2} f\big)^2
\cdot  P_0(\bar f)\notag\\
 &\quad +\frac12 \fe^{-it_n\partial_x^2}  P_N
 \Big[\big(\fe^{-it_n\partial_x^2}\bar f\big)
  \cdot \fe^{-i\tau\partial_x^2}
  P_N \big(\fe^{it_{n+1}\partial_x^2}
  \partial_x^{-1} f\big)^2\Big]\notag\\
&\quad -\frac12 \fe^{-it_n\partial_x^2}  P_N
 \Big[\big(\fe^{-it_n\partial_x^2}\bar f\big)
  \cdot
  P_N \big(\fe^{it_n\partial_x^2}
  \partial_x^{-1} f\big)^2\Big]\notag\\
  &\quad +i\tau\fe^{-it_n\partial_x^2}  P_N
 \Big[\big(\fe^{-it_n\partial_x^2}\bar f\big)
  \cdot
  P_N \big(\fe^{it_n\partial_x^2} f\big)^2\Big]\notag\\
 &\quad -2i\tau\fe^{-it_n\partial_x^2}  P_N
 \Big[\big(\fe^{-it_n\partial_x^2}\bar f\big)
  \cdot
  P_N \big(\fe^{it_n\partial_x^2} f\big)\Big]
 \cdot  P_0(f)\notag\\
& \quad +i\tau\fe^{-2it_n\partial_x^2}  P_N {\bar f}
 \cdot \big( P_0(f))^2.
\end{align}
Accordingly, for given $v^n\in S_N$, we compute $v^{n+1}\in S_N$ by
\begin{align}\label{NuSo-NLS}
v^{n+1}=\Phi_{\tau,N}^n\big(v^n\big),\qquad n\geq 0;\quad v^0=u^0.
\end{align}
This finishes the construction of the numerical scheme \eqref{numerical}.

\subsection{Some technical estimates}\label{subsec3}
We will frequently apply the following inequality, see \cite{Kato-Ponce}.
\begin{lemma}\label{lem:kato-Ponce}(Kato-Ponce inequality, \cite{Kato-Ponce}) The following estimates hold:
\begin{itemize}
  \item[(i)]
Let $f,g\in H^{\gamma}$ for some  $ \gamma>\frac 12$. Then we have
\begin{align*}
\|J^\gamma (fg)\|_{L^2}\lesssim \|f\|_{H^\gamma}\|g\|_{H^{\gamma}}.
\end{align*}
  \item[(ii)]
Let $f\in H^{\gamma+\gamma_1}$, $g\in H^{\gamma}$ for some $\gamma\ge 0$, $\gamma_1>\frac 12$. Then we have
\begin{align*}
\|J^\gamma (fg)\|_{L^2}\lesssim \|f\|_{H^{\gamma+\gamma_1}}\|g\|_{H^{\gamma}}.
\end{align*}
\end{itemize}
\end{lemma}

Next we present two specific estimates, which are used in section 4.
\begin{lemma}\label{lem:An}
The following bounds hold:
\begin{itemize}
\item[(i)]
Let $v\in {L^\infty((0,T);H^{\gamma})}$ for some $\gamma>\frac{1}{2}$, and $g\in T_3(1_{|k|> N};v)$. Then
$$
\| g\|_{L^2}
\lesssim N^{-\gamma}\|v\|_{L^\infty((0,T);H^{\gamma})}^3.
$$
\item[(ii)]
Let  $v\in {L^\infty((0,T);H^{\gamma})}$ for some $\gamma>\frac{1}{2}$, and
$g\in T_3(1_{|k_2+k_3|> N};v)$. Then
$$
\| g\|_{L^2}
\lesssim N^{-\gamma}\|v\|_{L^\infty((0,T);H^{\gamma})}^3.
$$
  \item[(iii)]
Let  $v\in {L^\infty((0,T);H^{\gamma})}$ for some $\gamma>\frac{1}{2}$, and
$g\in T_m(1;v)$, $m\geq1$. Then
$$
\| g\|_{H^{\gamma}}
\lesssim \|v\|_{L^\infty((0,T);H^{\gamma})}^m.
$$
\end{itemize}
\end{lemma}
\begin{proof}
We employ the notation $\hat V_{k_j}=|\hat v_{k_j}(t)|$.

(i) If $g\in T_3(1_{|k|> N};v)$,
using the definition of $ T_m(M;v)$ in \eqref{def:TmM},
 we have
\begin{align}\label{formula1}
|\hat g_k|
\lesssim
\left\{\begin{aligned}
&\sum\limits_{\substack{k=k_1+k_2+k_3}}
\>\hat{V}_{ k _1} \hat{V}_{ k _2}  \hat V_{ k _3}& \quad {\rm for} \quad |k|> N,\\
 &0 & \quad {\rm for} \quad |k|\leq N.
\end{aligned}\right.
\end{align}
By Parseval's identity,  we get
\begin{align*}
\|g\|_{L^2}^2=2\pi\sum_{k\in\Z}|\hat g_k|^2.
\end{align*}
From the above identity and \eqref{formula1}, the following
estimates hold
\begin{align*}
\|g\|_{L^2}&\lesssim
\bigg(\sum_{|k|> N} \Big(\sum\limits_{k=k_1+k_2+k_3}
\>\hat{V}_{ k _1} \hat{V}_{ k _2}  \hat V_{ k _3}\Big)^2\bigg)^\frac12\\
& \lesssim  N^{-\gamma}
\bigg(\sum_{|k|> N} \Big(\sum\limits_{k=k_1+k_2+k_3}|k|^\gamma
\>\hat{V}_{ k _1} \hat{V}_{ k _2}  \hat V_{ k _3}\Big)^2\bigg)^\frac12,
\end{align*}
where the last estimate holds true  because the frequency is limited to
$|k|>N$.

For convenience,  we employ the notation
$$
\tilde V=\sum_{k\in \Z} e^{ i  kx } \hat V_k.
$$
Then,  we have that $\widehat {\tilde V}_k=\hat  V_k=|\hat{v}_k(t)|$  and thus
\begin{align}\label{formu:V}
\|\tilde V\|_{H^\gamma}^2
=&2\pi\sum_{k\in \Z}(1+ k ^2)^\gamma |\widehat {\tilde V}_k|^2
=2\pi\sum_{k\in \Z}(1+ k ^2)^\gamma |\hat  V_k|^2\notag\\
=& 2\pi\sum_{k\in \Z}(1+ k ^2)^\gamma |\hat{v}_k(t)|^2
=\|v\|_{H^\gamma}^2.
\end{align}
Using Parseval's identity once more we obtain
\begin{align*}
\|g\|_{L^2} & \lesssim N^{-\gamma}
  \big\| \tilde V^3\big\|_{L^\infty((0,T);H^\gamma)}
  \lesssim N^{-\gamma}\|\tilde V\|_{L^\infty((0,T);H^{\gamma})}^3.
\end{align*}
Therefore, from   \eqref{formu:V} and the above inequality, we have
\begin{align*}
 \|g\|_{L^2}  \lesssim N^{-\gamma}\|v\|_{L^\infty((0,T);H^{\gamma})}^3.
\end{align*}

(ii) We use the same argument as in (i) and Parseval's identity to get
\begin{align*}
\| g\|_{L^2}
&\lesssim
\bigg(\sum_{k\in\Z} \Big(\sum\limits_{\substack{k=k_1+k_2+k_3\\|k_2+k_3|>N}}
\>\hat{V}_{ k _1} \hat{V}_{ k _2}  \hat V_{ k _3}\Big)^2\bigg)^\frac12\\
 &\lesssim  \big\| \tilde V\cdot P_{>N}(\tilde V^2)\big\|_{L^\infty((0,T);L^2)},
\end{align*}
where the operator $P_{>N}$ is defined in \eqref{truncation}.

Then, by H\"{o}lder's inequality, we have
\begin{align*}
\| g\|_{L^2}& \lesssim
  \big\|\tilde V\big\|_{L^\infty((0,T);L^\infty)}
  \big\| P_{>N}(\tilde V^2)\big\|_{L^\infty((0,T);L^2)}.
  \end{align*}
  Note that  $$P_{>N}(\tilde V^2)\in T_2(1_{|k|> N};v),$$
  hence   we use the result of (i) and \eqref{formu:V} to obtain
  $$
   \big\| P_{>N}(\tilde V^2)\big\|_{L^\infty((0,T);L^2)}
   \lesssim N^{-\gamma}
  \big\| v\big\|_{L^\infty((0,T);H^\gamma)}^2.
  $$
Now,  we can use Sobolev's inequality to finally obtain
  \begin{align*}
  \| g\|_{L^2}\lesssim N^{-\gamma}
  \big\| v\big\|_{L^\infty((0,T);H^\gamma)}^3.
\end{align*}

(iii)Using the explicit form of  $ T_m(M;v)$ given in \eqref{def:TmM}, we have
\begin{align*}
|\hat g_k|\lesssim
\sum\limits_{ k = k _1+\cdots+ k _m}
\>\hat{V}_{ k _1} \cdots \hat V_{ k _m}.
\end{align*}
By Parseval's identity, Lemma \ref{lem:kato-Ponce} and \eqref{formu:V} again, we obtain that
\begin{align*}
\big\| g\big\|_{H^{\gamma}}
  \lesssim \|\tilde V^m\|_{L^\infty((0,T);H^{\gamma})}
\lesssim \|v\|_{L^\infty((0,T);H^{\gamma})}^m.
\end{align*}
This concludes the proof.
\end{proof}

Furthermore, from Lemma \ref{lem:An}, we have for any $\gamma>\frac12$,
\begin{align}\label{r1}
&\big\|\mathcal{R}_1(v)\big\|_{L^2}
\leq\big\|\mathcal{R}_1(v)\big\|_{H^{\gamma}}
\lesssim \tau^2\|v\|_{L^\infty((0,T);H^{\gamma})}^5,\\
&\big\|\mathcal{R}_2(v)\big\|_{L^2}
\lesssim \tau N^{-\gamma} \|v\|_{L^\infty((0,T);H^{\gamma})}^3,\qquad
\big\|\mathcal{R}_2(v)\big\|_{H^{\gamma}}
\lesssim \tau\|v\|_{L^\infty((0,T);H^{\gamma})}^3,\label{r2}\\
&\big\|\mathcal{R}_4(v)\big\|_{L^2}
\lesssim \tau N^{-\gamma} \|v\|_{L^\infty((0,T);H^{\gamma})}^3,\qquad
\big\|\mathcal{R}_4(v)\big\|_{H^{\gamma}}
\lesssim \tau\|v\|_{L^\infty((0,T);H^{\gamma})}^3\label{r4}.
\end{align}

\begin{lemma}\label{lem:r3}
The following estimates hold:
\begin{itemize}
\item[(i)]
Let $v\in {L^\infty((0,T);H^{\gamma})}$ for some $1\geq\gamma>\frac{1}{2}$. Then for any small enough $\varepsilon>0$,
$$
\big\|\mathcal{R}_3(v)\big\|_{L^2}
\lesssim \tau^{\frac32\gamma+\frac12-\varepsilon}\|v\|_{L^\infty((0,T);H^{\gamma})}^3.
$$
  \item[(ii)]
Let $v\in {L^\infty((0,T);H^{\gamma})}$ for some $\gamma>\frac{1}{2}$. Then,
$$
\big\|\mathcal{R}_3(v)\big\|_{H^{\gamma}}
\lesssim \tau\|v\|_{L^\infty((0,T);H^{\gamma})}^3.
$$
\end{itemize}
\end{lemma}
\begin{proof}
(i) As employed in Lemma 3.2, we use the notation $\hat V_{k_j}=|\hat v_{k_j}(t)|$.
Note that for $s\in [0,\tau]$,
$$
 \big|\fe^{2iskk_1}-1\big|
 \lesssim \tau^\frac{\gamma}{2}|k|^\frac{\gamma}{2}|k_1|^\frac{\gamma}{2},
 \qquad  \big|\fe^{2isk_2k_3}-1\big|
 \lesssim \tau^{\gamma-\frac12-\varepsilon}
 |k_2|^{\gamma-\frac12-\varepsilon}|k_3|^{\gamma-\frac12-\varepsilon}.
$$
Therefore, we obtain from \eqref{def:r3} that
\begin{align*}
\big|\hat{ \mathcal{R}}_{3,k}(v)\big|
\lesssim
 \tau^{\frac32\gamma+\frac12-\varepsilon}
 \sum\limits_{{\substack{k=k_1+k_2+k_3\\   |k_2+k_3|\leq N}}}
|k|^\frac{\gamma}{2}|k_1|^\frac{\gamma}{2}|k_2|^{\gamma-\frac12-\varepsilon}|k_3|^{\gamma-\frac12-\varepsilon}
      \>\hat{V}_{ k _1} \hat{V}_{ k _2}  \hat V_{ k _3}.
\end{align*}
Based on the relation between the frequencies, we consider three cases.

Case 1: $|k|\lesssim|k_1|$.
From the above estimate, we have
\begin{align*}
\big|\hat{ \mathcal{R}}_{3,k}(v)\big|
\lesssim
 \tau^{\frac32\gamma+\frac12-\varepsilon}
 \sum\limits_{{\substack{k=k_1+k_2+k_3\\   |k_2+k_3|\leq N}}}
|k_1|^\gamma|k_2|^{\gamma-\frac12-\varepsilon}|k_3|^{\gamma-\frac12-\varepsilon}
     \>\hat{V}_{ k _1} \hat{V}_{ k _2}  \hat V_{ k _3}.
\end{align*}
We denote as before
$$
\tilde V=\sum_{k\in \Z} e^{ i  kx } \hat V_k.
$$
Then from \eqref{formu:V}, we  have
$
\|\tilde V\|_{H^\gamma}
=\|v\|_{H^\gamma}.
$
Therefore, by Parseval's identity, we obtain that, for any $\gamma>\frac12$,
\begin{align*}
\big\|\mathcal{R}_3(v)\big\|_{L^2}
\lesssim& \tau^{\frac32\gamma+\frac12-\varepsilon}
\Big\||\partial_x|^\gamma\tilde V
\cdot \big(|\partial_x|^{\gamma-\frac12-\varepsilon}\tilde V\big)^2\Big\|_{L^2}\\
\lesssim &
\tau^{\frac32\gamma+\frac12-\varepsilon}
\|v\|_{L^\infty((0,T);H^{\gamma})}^3.
\end{align*}

Case 2: $|k|\lesssim|k_2|$. In this case, we have
\begin{align*}
\big|\hat{ \mathcal{R}}_{3,k}(v)\big|
\lesssim
 \tau^{\frac32\gamma+\frac12-\varepsilon}
 \sum\limits_{{\substack{k=k_1+k_2+k_3\\   |k_2+k_3|\leq N}}}
|k_1|^\frac{\gamma}{2}
|k_2|^{\frac32\gamma-\frac12-\varepsilon}|k_3|^{\gamma-\frac12-\varepsilon}
       \>\hat{V}_{ k _1} \hat{V}_{ k _2}  \hat V_{ k _3}.
\end{align*}
Therefore, we use Parseval's identity to obtain
\begin{align*}
\big\|\mathcal{R}_3(v)\big\|_{L^2}
\lesssim \tau^{\frac32\gamma+\frac12-\varepsilon}
\Big\||\partial_x|^{\frac12\gamma}\tilde V
\cdot|\partial_x|^{\frac32\gamma-\frac12-\varepsilon}\tilde V
\cdot|\partial_x|^{\gamma-\frac12-\varepsilon}\tilde V\Big\|_{L^2}.
\end{align*}
If $\gamma\neq1$, we employ the H\"{o}lder and Sobolev inequalities to get
\begin{align*}
\big\|\mathcal{R}_3(v)\big\|_{L^2}
&\lesssim \tau^{\frac32\gamma+\frac12-\varepsilon}
\big\||\partial_x|^{\frac12\gamma}\tilde V\big\|_{L^{\frac{2}{1-\gamma}}}
\cdot\big\||\partial_x|^{\frac32\gamma-\frac12-\varepsilon}\tilde V\big\|_{L^\frac2\gamma}
\cdot\big\||\partial_x|^{\gamma-\frac12-\varepsilon}\tilde V\big\|_{L^\infty}\\
&\lesssim \tau^{\frac32\gamma+\frac12-\varepsilon}
\|v\|_{L^\infty((0,T);H^{\gamma})}^3.
\end{align*}
If $\gamma=1$,  we write
\begin{align*}
\big\|\mathcal{R}_3(v)\big\|_{L^2}
&\lesssim \tau^{2-\varepsilon}
\big\||\partial_x^{\frac12}| \tilde V\big\|_{L^{\frac{1}{\varepsilon}}}
\cdot\big\||\partial_x^{1-\varepsilon}|\tilde V\big\|_{L^\frac{2}{1-2\varepsilon}}
\cdot\big\||\partial_x^{\gamma-\frac12-\varepsilon}|\tilde V\big\|_{L^\infty}\\
&\lesssim \tau^{2-\varepsilon}\|v\|_{L^\infty((0,T);H^{\gamma})}^3.
\end{align*}
Hence, for any $\frac12<\gamma\leq1$ we conclude that
$$
\big\|\mathcal{R}_3(v)\big\|_{L^2}
\lesssim \tau^{\frac32\gamma+\frac12-\varepsilon}
\|v\|_{L^\infty((0,T);H^{\gamma})}^3.
$$

Case 3: $|k|\lesssim|k_3|$. It is same as case 2.

(ii) We have the following inequality
$$
 \Big|\big(\fe^{2iskk_1}-1\big)\big(\fe^{2isk_2k_3}-1\big)\Big|
 \lesssim1.
$$
 Plugging the above inequality into \eqref{def:r3}, we obtain
\begin{align}\label{for3}
\big|\hat{ \mathcal{R}}_{3,k}(v)\big|
\lesssim
 \tau
 \sum\limits_{{\substack{k=k_1+k_2+k_3\\   |k_2+k_3|\leq N}}}
        \>\hat{V}_{ k _1} \hat{V}_{ k _2}  \hat V_{ k _3}.
\end{align}
Then, from the definition of $ T_m(M;v)$ in \eqref{def:TmM}, we  observe   that
$$
\mathcal{R}_{3}(v)
\in\tau T_3(1;v).
$$
Finally, from Lemma \ref{lem:An} (iii), we get the result  (ii).
\end{proof}

\section{Proof of Theorem \ref{result1}}\label{proof}

Note that the $L^2$ stability estimate depends on bounds
of the numerical solution in $H^\gamma$, $\frac12<\gamma\leq1$.
Therefore we first show the $H^\gamma$  bound before we give the proof of Theorem \ref{result1}.
 This strategy was first used in  \cite{Lubich}.
In particular,  the bound of the $H^\gamma$ norm of the numerical solution
is independent of the degrees of freedom in the spatial discretisation, $N$. It only depends  on $T$ and the bound $\|u\|_{L^\infty((0,T);H^{\gamma})}$.
We are now in a position to show the bound in $H^\gamma$.

\begin{proposition}\label{bound}
Let $u^n_{\tau,N}$ be the numerical solution given in \eqref{numerical}. Under the assumption that $u^0\in H^{\gamma}(\bT)$ for some $\frac{1}{2}<\gamma\leq1$,  there exists a positive constant $C$, such that for any
integer $N\geq1$
\begin{equation}
  \|u^n_{\tau,N}\|_{H^\gamma}\leq C, \quad {\text{for all}} \quad 0\leq n\tau\leq T,
\end{equation}
 where the constant $C$ only depends  on $T$ and the bound $\|u\|_{L^\infty((0,T);H^{\gamma})}$.
\end{proposition}
\begin{proof}
Let $v^n=u^n_{\tau,N}$. From \eqref{NuSo-NLS}, we have
\begin{align}\label{err}
v(t_{n+1})-v^{n+1}&=v(t_{n+1})-\Phi^n_{\tau,N}\big(v(t_n)\big)+ \Phi^n_{\tau,N}\big(v(t_n)\big)-\Phi^n_{\tau,N}\big(v^n\big)\\
          &=\mathcal{L}^n+\Phi^n_{\tau,N}\big(v(t_n)\big)-\Phi^n_{\tau,N}\big(v^n\big),
\end{align}
where $\mathcal{L}^n=v(t_{n+1})-\Phi^n_{\tau,N}\big(v(t_n)\big)$.

Furthermore, from \eqref{vtn+1-app} we get
$$
\mathcal{L}^n=\mathcal{R}_1(v)+\mathcal{R}_2(v)
+\mathcal{R}_3(v)+\mathcal{R}_4(v).
$$
Then, from \eqref{r1}, \eqref{r2}, \eqref{r4} and
Lemma \ref{lem:r3} (ii), we have
\begin{align}\label{est:ln}
\|\mathcal{L}^n\|_{H^\gamma}\leq C\tau.
\end{align}
Note that the constant  $C$ only depends   on $\|u\|_{L^\infty((0,T);H^{\gamma})}$.

Recall that   $\Phi^n_{\tau,N}(f)$ defined in \eqref{Phi-def} can be written the following integral form:
\begin{align}\label{stab}
\widehat\Phi^n_{\tau,N}(f)(k)=\hat f_k-i\sum\limits_{{\substack{k=k_1+k_2+k_3\\   |k_2+k_3|\leq N\\|k|\leq N}}}
\int_0^\tau \fe^{it_n\phi}\Big[\fe^{2iskk_1}+\big(\fe^{2isk_2k_3}-1\big)\Big]\,ds
       \>\widehat{ \bar f}_{k_1}\hat f_{k_2}\hat f_{k_3}.
\end{align}
Further, note that
\begin{align*}
\Big|\fe^{it_n\phi}\Big[\fe^{2iskk_1}+\big(\fe^{2isk_2k_3}-1\big)\Big]\Big|\lesssim1.
\end{align*}
Then, by Lemma \ref{lem:kato-Ponce} (i), we obtain
\begin{align}\label{stability3}
\|\Phi^n_{\tau,N}\big(v(t_n)\big)-\Phi^n_{\tau,N}\big(v^n\big)\|_{H^\gamma}
   \leq &(1+C\tau) \|v^n-v(t_n)\|_{H^\gamma}
        +C\tau\|v^n-v(t_n)\|_{H^\gamma}^3.
   \end{align}
A combination of  the above estimates allows that
\begin{align*}
\|v(t_{n+1})-v^{n+1}\|_{H^\gamma}
   \leq (1+C\tau) \|v^n-v(t_n)\|_{H^\gamma}
        +C\tau\|v^n-v(t_n)\|_{H^\gamma}^3+C\tau.
\end{align*}
By iteration and Gronwall's lemma, we get
\begin{align*}
\big\|v(t_{n+1})-v^{n+1}\big\|_{H^\gamma}
\le C.
\end{align*}
That means
\begin{align*}
\big\|v^{n+1}\big\|_{H^\gamma}
\le C.
\end{align*}
This finishes the proof on the boundness in $H^\gamma$.
\end{proof}

Now we start to prove  Theorem \ref{result1}.

From \eqref{r1}, \eqref{r2}, \eqref{r4} and
Lemma \ref{lem:r3} (i), we have
\begin{align}\label{est:ln2}
\|\mathcal{L}^n\|_{L^2}\leq C\tau^{\frac32\gamma+\frac12-\varepsilon}
+C\tau N^{-\gamma},
\end{align}
where the constant  $C$ only depends   on the bound $\|u\|_{L^\infty((0,T);H^{\gamma})}$.

From \eqref{stab},
we have
\begin{align*}
&\|\Phi^n_{\tau,N}\big(v(t_n)\big)-\Phi^n_{\tau,N}\big(v^n\big)\|_{L^2}\\
  &\qquad\qquad \leq \|v^n-v(t_n)\|_{L^2}
   +C\tau\|v^n-v(t_n)\|_{L^2}
   \big(\|v^n\|_{H^\gamma}^2+\|v(t_n)\|_{H^\gamma}^2\big).
   \end{align*}
By Proposition \ref{bound}, we have
\begin{align*}
\|\Phi^n_{\tau,N}\big(v(t_n)\big)-\Phi^n_{\tau,N}\big(v^n\big)\|_{L^2}
   \leq
   (1+C\tau)\|v^n-v(t_n)\|_{L^2}.
   \end{align*}
Combining \eqref{err} with \eqref{est:ln2} and the above estimate, we get
\begin{align*}
\|v(t_{n+1})-v^{n+1}\|_{H^\gamma}
   \leq &(1+C\tau) \|v^n-v(t_n)\|_{H^\gamma}
        +C\tau^{\frac32\gamma+\frac12-\varepsilon}
+C\tau N^{-\gamma}.
\end{align*}
By iteration and Gronwall's lemma, we finally get
\begin{align*}
\big\|v(t_{n+1})-v^{n+1}\big\|_{H^\gamma}
\le C\tau^{\frac32\gamma-\frac12-\varepsilon}
+C N^{-\gamma}.
\end{align*}
This concludes the proof. 

\section{Numerical experiments} \label{sec:numerical}
In this section we carry out numerical experiments to support our theoretical analysis.
We consider the
nonlinear Schr\"{o}dinger equation \eqref{model} with initial data
\begin{align}\label{initial}
u^0(x)=\sum_{k\in\Z}(1+|k|)^{-\frac12-\gamma}g_k \fe^{ikx},
\end{align}
where $\gamma$ and $(g_k)_{k\in\Z}$ are used to set the regularity of the data. The complex coefficients $g_k$ are chosen as uniformly
distributed random variables in $[-1,1]+i[-1,1]$.
They are generated with the matlab routine rand.
This choice guarantees that $u^0\in H^\gamma$.
however, as a consequence of the Paley-Zygmund Theorem
(see, e.g., \cite{kahane}), the initial data \eqref{initial}
satisfies the stronger regularity condition
$\partial_x^\gamma u^0\in L^p$, $2\leq p< \infty$.
This will be a slight issue when we study the temporal discretisation
error. See the discussion below.

We start our numerical experiments with the spatial  discretisation errors. On the left-hand side of Fig.~1, we present our results
for a sufficiently small time step size $\tau$. This allows us to ignore the errors caused by temporal discretisation.
We choose three different values of $\gamma\in [\frac12,1]$ to illustrate the spatial convergence rate of our scheme.
In order to measure the spatial discretization error
$u(t_m,\cdot)-u^m_{\tau,N}$ for fixed time $t_m$, we use
 the discrete $L^2$ norm
\begin{align}\label{l2}
\|w\|_{L^2_N}^2=\frac{N}{2\pi}\sum_{j=0}^{N-1}|w(x_j)|^2,
\qquad x_j=\frac{2j}{N}\pi.
\end{align}
The results of our numerical experiments agree well with the corresponding results of the theoretical analysis,
see Theorem \ref{result1}.

The temporal discretisation errors are displayed on the right-hand side of Fig.~ 1 for a sufficiently large $N$.
This allows us to ignore the errors caused by spatial discretisation.
In order to illustrate the time convergence rate, we present results for three different values of $\gamma\in [\frac12,1]$.
We use again the norm specified in \eqref{l2}.
Note that the observed rates of convergence are slightly better
than predicted by Theorem \ref{result1}.
The reasons for this is the additional regularity,
guaranteed by the  Paley-Zygmund Theorem  \cite{kahane}.

\begin{figure}[!htb]
$$\begin{array}{cc}
\psfig{figure=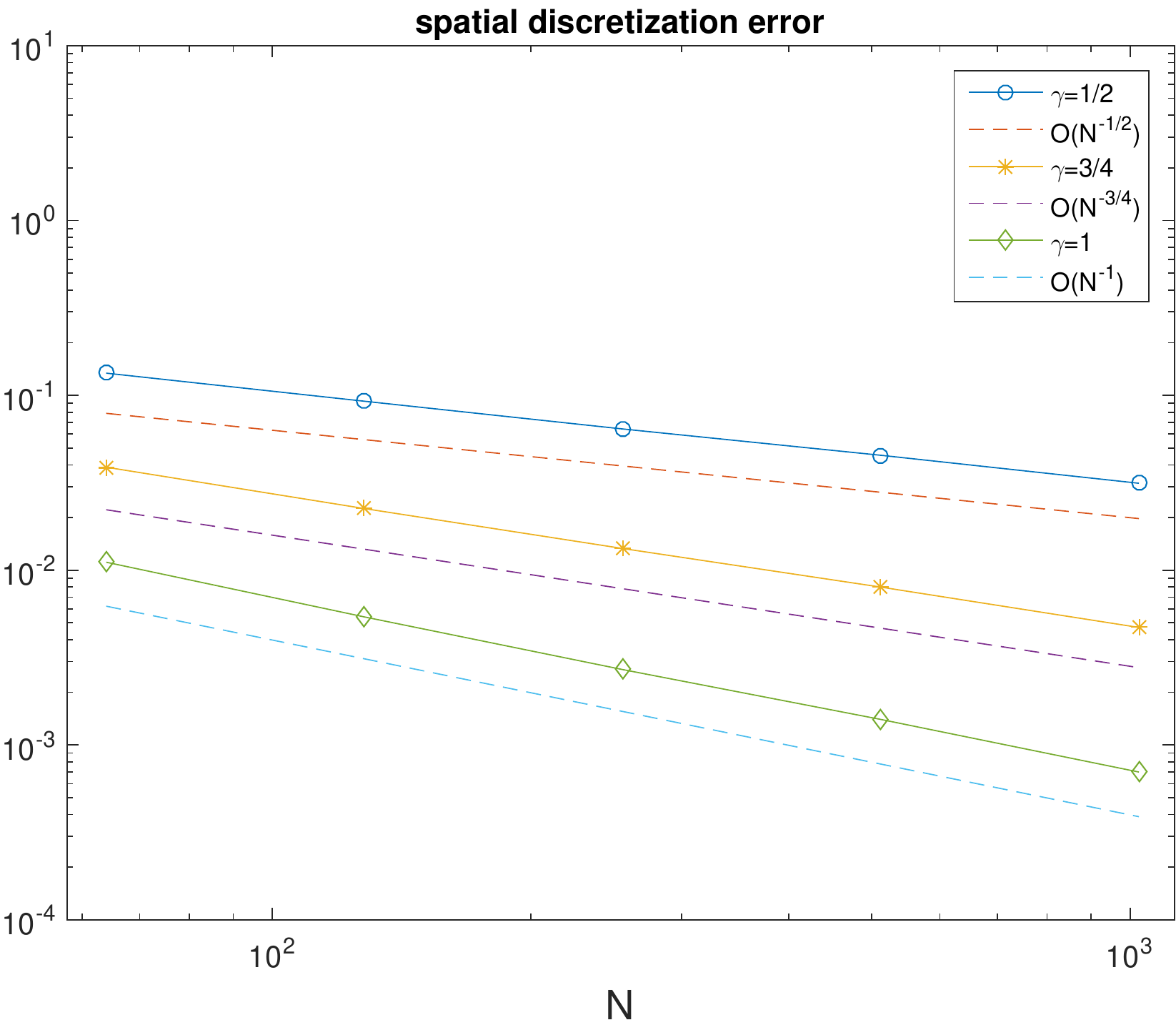,height=7.2cm,width=6.8cm}\hspace{12mm}
\psfig{figure=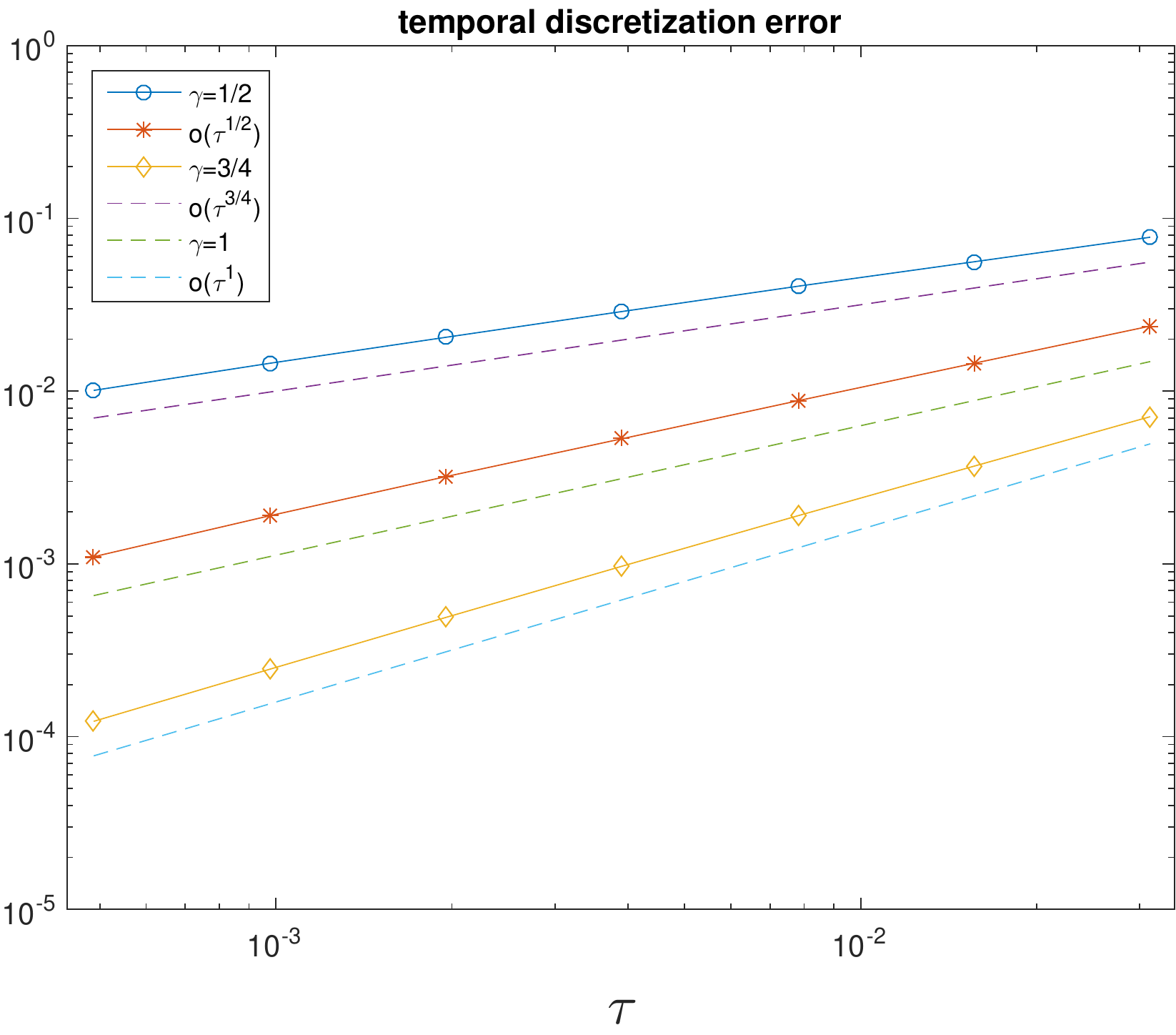,height=7.2cm,width=6.8cm}
\end{array}$$
\caption{Spatial discretisation error  at $\tau=2^{-15}$, $T=1$
	for various values of $N$ and $\gamma$ (left)
and temporal discretisation error  for $N=2^{14}$ at $T=1$
for various values of $\tau$ and $\gamma$ (right).}
\label{fig:NLRI}
\end{figure}

\section{Conclusion} \label{sec:conclusion}
 We have constructed a fully discrete low-regularity integrator
for  the cubic NLS equation with nonsmooth initial data in one space dimension.
The scheme can be computed with FFT   with
$\mathcal{O}(N\log N)$ operations  per time step.
We have proved convergence  in $L^2(\T)$ for
initial data in $H^{\gamma}(\T)$, $\frac12<\gamma\leq 1$.
Numerical results illustrate our convergence result.

\section*{acknowledgements}
	\noindent This research is supported by the NSFC key project under the grant number 11831003, and NSFC under the grant numbers 11971356.
	The second author also acknowledges financial support by the
	China Scholarship Council.


\vskip 25pt
\bibliographystyle{model1-num-names}

\begin{thebibliography}{00}

\bibitem{besse}
{\sc C. Besse, B. Bid\'{e}garay, and S. Descombes},
Order estimates in time of splitting methods
for the nonlinear Schr\"odinger equation, SIAM J. Numer. Anal. 40 (2002), pp. 26--40.


\bibitem{Bo}
{\sc J. Bourgain}, Fourier transform restriction phenomena for certain lattice subsets and
applications to nonlinear evolution equations. I. Schr\"odinger  equations, Geom. Funct. Anal.
 3  (1993), pp. 107--156.




%
%
%

\bibitem{faou}
{\sc E. Faou,}
 Geometric Numerical Integration and Schr\"odinger  Equations,
  European Mathematical Society Publishing House, Z\"{u}rich, 2012.




\bibitem{jahnke}
{\sc T. Jahnke and C. Lubich,}
 Error bounds for exponential operator splittings, BIT 40 (2000), pp. 735--744.

 \bibitem{kahane}
 {\sc J. P. Kahane,}
 Some random series of functions, Cambridge University Press, Cambridge, 1993.

\bibitem{Kato-Ponce}
{\sc T. Kato and G. Ponce}, Commutator estimates and the Euler and Navier-Stokes equations, Commun. Pure Appl. Math. 41 (1988), pp. 891-907.




\bibitem{lownls2}
{\sc M. Kn\"{o}ller, A. Ostermann, and K. Schratz},
A Fourier integrator for the cubic nonlinear Schr\"{o}dinger equation with rough initial data,
SIAM J. Numer. Anal. 57 (2019), pp. 1967--1986.


\bibitem{liwu}
{\sc B. Li and Y. Wu},
A full discrete low-regularity integrator for the 1D period cubic nonlinear Schr\"{o}dinger equation,
arXiv:2101.03728, 2021.

\bibitem{Lubich}
{\sc C. Lubich}, On splitting methods for Schr\"{o}dinger-Poisson and cubic nonlinear Schr\"{o}dinger equations, Math. Comp. 77 (2008), pp. 2141--2153.

\bibitem{lownls} {\sc A. Ostermann and K. Schratz},
{Low regularity exponential-type integrators for semilinear Schr\"{o}dinger equations},
 Found. Comput. Math. 18 (2018), pp. 731--755.

 \bibitem{ostermann}
{\sc A. Ostermann, F. Rousset, and K. Schratz},
Error estimates of a Fourier integrator for the cubic Schr\"{o}dinger equation at low regularity,
 Found. Comput. Math. 21 (2021), pp. 725--765.


\bibitem{ostermann2020}
{\sc A. Ostermann, F. Rousset, and K. Schratz},
Fourier integrator for periodic NLS: low regularity estimates via discrete Bourgain spaces,
 arXiv:2006.12785, 2020.

\bibitem{sanz}
{\sc J. M. Sanz-Serna,}
Methods for the numerical solution of the nonlinear Schr\"{o}dinger equation, Math. Comp. 43 (1984), pp. 21--27.



\bibitem{w}
{\sc Y. Wu and F. Yao,}
A first-order Fourier integrator for the nonlinear Schr\"odinger equation on $\T$ without loss of regularity,
arXiv:2010.02672, 2020.


\end{thebibliography}

\end{document}